\documentclass[11pt]{article}
\usepackage{amsmath,amsthm,amsfonts}
\newcommand\pd[2]{\frac{\partial#1}{\partial#2}}
\renewcommand{\=}{\doteq}
\newcommand{\p}{\partial}
\newcommand{\X}{\mathfrak{X}}

\newcommand{\la}{\lambda}
\newtheorem{thm}{Theorem}[section]

 \newtheorem{lemma}[thm]{Lemma}
\theoremstyle{definition}

\theoremstyle{definition}
 
\numberwithin{equation}{section}
\numberwithin{equation}{section}
\begin{document}
\title{\bf  Moufang symmetry XI.\\
Integrability of generalized Lie equations of\\ continuous Moufang transformations
}
\author{Eugen Paal}
\date{}
\maketitle
\thispagestyle{empty}
\begin{abstract}
Integrability of generalized Lie equations of continuous Moufang transformations is inquired. 
\end{abstract}

\section{Introduction}

In this paper we proceed explaing the Moufang symmetry. The present paper can be seen as a continuation of  \cite{Paal1,Paal2,Paal_AAM}. 

\section{Generalized Lie equations}

In \cite{Paal_AAM} the \emph{generalized Lie equations} (GLE)  of continous Moufang transformations were found. For $S_gA$ the GLE read
\begin{subequations}
\label{gle_S}
\begin{align}
u^{s}_{j}(g)\pd{(S_gA)^{\mu}}{g^{s}}+T^{\nu}_{j}(A)\pd{(S_gA)^{i}}{A^{\nu}}+P^{\nu}_{j}(S_gA)&=0\\
v^{s}_{j}(g)\pd{(S_gA)^{\mu}}{g^{s}}+P^{\nu}_{j}(h)\pd{(S_gA)^{i}}{A^{\nu}}+T^{\nu}_{j}(S_gA)&=0\\
w^{s}_{j}(g)\pd{(S_gA)^{\mu}}{g^{s}}+S^{\nu}_{j}(h)\pd{(S_gA)^{i}}{A^{\nu}}+S^{\nu}_{j}(S_gA)&=0
\end{align}
\end{subequations}
where $gh$ is the product of $g$ and $h$, and the auxiliary functions $u^s_j$, $v^s_j$, $w^s_j$ and
$S^\mu_j$, $T^\mu_j$, $P^\mu_j(g)$ are related with the  constraints
\begin{gather}
u^s_j(g)+v^s_j(g)+w^s_j(g)=0\\
S^\mu_j(A)+T^\mu_j(A)+P^\mu_j(A)=0
\end{gather}
For $T_gA$ the GLE read
\begin{subequations}
\label{gle_T}
\begin{align}
v^{s}_{j}(g)\pd{(T_gA)^{\mu}}{g^{s}}+S^{\nu}_{j}(A)\pd{(T_gA)^{i}}{A^{\nu}}+P^{\nu}_{j}(T_gA)&=0\\
u^{s}_{j}(g)\pd{(T_gA)^{\mu}}{g^{s}}+P^{\nu}_{j}(h)\pd{(T_gA)^{i}}{A^{\nu}}+S^{\nu}_{j}(T_gA)&=0\\
w^{s}_{j}(g)\pd{(T_gA)^{\mu}}{g^{s}}+T^{\nu}_{j}(h)\pd{(T_gA)^{i}}{A^{\nu}}+T^{\nu}_{j}(T_gA)&=0
\end{align}
\end{subequations}
In this paper we inquire integrability of GLE (\ref{gle_S}a--c) and (\ref{gle_T}a--c).

\section{Generalized Maurer-Cartan equations and Yamagutian I}

Recall from \cite{Paal1} that for $x$ in  $T_e(G)$ the infinitesimal translations of $G$ are defined by
\begin{equation*}
L_x\=x^j u^s_j(g)\pd{}{g^s},\quad
R_x\=x^j v^s_j(g)\pd{}{g^s},\quad
M_x\=x^j w^s_j(g)\pd{}{g^s}\quad \in T_g(G)
\end{equation*}
with constriant
\begin{equation*}
L_x+R_x+M_x=0
\end{equation*}
Following triality \cite{Paal2} define the Yamagutian $Y(x;y)$ by
\begin{equation*}
6Y(x;y)=[L_x,L_y]+[R_x,R_y]+[M_x,M_y]
\end{equation*}
We know  from \cite{Paal2} the generalized Maurer-Cartan equations:
\begin{subequations}
\label{m-c}
\begin{align} 
[L_{x},L_{y}]&=L_{[x,y]}-2[L_{x},R_{y}]\\
[R_{x},R_{y}]&=R_{[y,x]}-2[R_{x},L_{y}]\\
[L_{x},R_{y}]&=[R_{x},L_{y}],\quad \forall x,y \in T_e(G) 
\end{align}
\end{subequations}
The latter can be written \cite{Paal2} as follows:
\begin{subequations}
\label{lr-yam}
\begin{align}
[L_{x},L_{y}]&=2Y(x;y)+\frac{1}{3}L_{[x,y]}+\frac{2}{3}R_{[x,y]}\\
[L_{x},R_{y}]&=-Y(x;y)+\frac{1}{3}L_{[x,y]}-\frac{1}{3}R_{[x,y]}\\
[R_{x},R_{y}]&=2Y(x;y)-\frac{2}{3}L_{[x,y]}-\frac{1}{3}R_{[x,y]}
\end{align}
\end{subequations}
Define the  (secondary) auxiliary functions of $G$ by
\begin{align*}
u^s_{jk}(g)
&\=u^p_k(g)\pd{u^s_j(g)}{g^p}-u^p_j(g)\pd{u^s_k(g)}{g^p}\\
v^s_{jk}(g)
&\=v^p_k(g)\pd{v^s_j(g)}{g^p}-v^p_j(g)\pd{v^s_k(g)}{g^p}\\
w^s_{jk}(g)
&\=w^p_k(g)\pd{w^s_j(g)}{g^p}-w^p_j(g)\pd{w^s_k(g)}{g^p}
\end{align*}
The Yamaguti functions $Y^i_{jk}$ are defined by
\begin{equation*}
6Y^s_{jk}(g)\=u^s_{jk}(g)+v^s_{jk}(g)+w^s_{jk}(g)
\end{equation*}
Evidently,
\begin{align*}
[L_x,L_y]
&=-x^jy^k u^s_{jk}(g)\pd{}{g^s}\\
[R_x,R_y]
&=-x^jy^k v^s_{jk}(g)\pd{}{g^s}\\
[M_x,M_y]
&=-x^jy^k w^s_{jk}(g)\pd{}{g^s}
\end{align*}
By adding the above formulae, we get
\begin{equation*}
Y(x;y)=-x^jy^k Y^s_{jk}(g)\pd{}{g^s}
\end{equation*}

\begin{lemma}
One has
\begin{subequations}
\label{uvw}
\begin{align}
u^i_{jk}
&\=2Y^i_{jk}+\frac{1}{3}C^s_{jk}(u^i_s+2v^i_s)\\
v^i_{jk}
&\=2Y^i_{jk}-\frac{1}{3}C^s_{jk}(2u^i_s+v^i_s)\\
w^i_{jk}
&\=2Y^i_{jk}+\frac{1}{3}C^s_{jk}(u^i_s-v^i_s)
\end{align}
\end{subequations}
\end{lemma}

\begin{proof}
To see (\ref{uvw}a,b) use (\ref{lr-yam}a,c) . To see (\ref{uvw}c) calculate by using (\ref{lr-yam}):
\begin{align*}
[M_x,M_y]
&=[L_x+R_x,L_y+R_y]\\
&=[L_x,L_y]+[L_x,R_y]+[R_x,L_y]+[R_x,R_y]\\
&=[L_x,L_y]+2[L_x,R_y]+[R_x,R_y]\\
&=2Y(x;y)+\frac{1}{3}\left(L_{[x,y]}-R_{[x,y]}\right)
\end{align*}
and  (\ref{uvw}b) easily follows.
\end{proof}

\section{Generalized Maurer-Cartan equations and Yamagutian II}

Recall from \cite{Paal_AAM} that for $x$ in  $T_e(G)$ the infinitesimal translations of $G$ are defined by
\begin{equation*}
S_x\=x^j S^\nu_j(A)\pd{}{A^\nu},\quad
T_x\=x^j T^\nu_j(A)\pd{}{A^\nu},\quad
P_x\=x^j P^\nu_j(A)\pd{}{A^\nu}\quad \in T_A(\X)
\end{equation*}
with constriant
\begin{equation*}
S_x+T_x+P_x=0
\end{equation*}
Following triality \cite{Paal2} define the Yamagutian $Y(x;y)$ by
\begin{equation*}
6Y(x;y)=[S_x,S_y]+[T_x,T_y]+[P_x,P_y]
\end{equation*}
We know  from \cite{Paal_AAM} the generalized Maurer-Cartan equations:
\begin{subequations}
\label{m-c_ST}
\begin{align} 
[S_{x},S_{y}]&=S_{[x,y]}-2[S_{x},T_{y}]\\
[T_{x},T_{y}]&=T_{[y,x]}-2[T_{x},S_{y}]\\
[S_{x},T_{y}]&=[T_{x},S_{y}],\quad \forall x,y \in T_e(G) 
\end{align}
\end{subequations}
The latter can be written \cite{Paal2} as follows:
\begin{subequations}
\label{lr-yam_ST}
\begin{align}
[S_{x},S_{y}]&=2Y(x;y)+\frac{1}{3}S_{[x,y]}+\frac{2}{3}T_{[x,y]}\\
[T_{x},T_{y}]&=-Y(x;y)+\frac{1}{3}S_{[x,y]}-\frac{1}{3}T_{[x,y]}\\
[T_{x},T_{y}]&=2Y(x;y)-\frac{2}{3}S_{[x,y]}-\frac{1}{3}T_{[x,y]}
\end{align}
\end{subequations}
Define the  (secondary) auxiliary functions of $G$ by
\begin{align*}
S^\mu_{jk}(A)
&\=S^\nu_k(A)\pd{S^\mu_j(A)}{A^\nu}-S^\nu_j(g)\pd{S^\mu_k(A)}{A^\nu}\\
T^\mu_{jk}(A)
&\=T^\nu_k(A)\pd{T^\mu_j(A)}{A^\nu}-T^\nu_j(g)\pd{T^\mu_k(A)}{A^\nu}\\
P^\mu_{jk}(A)
&\=P^\nu_k(A)\pd{P^\mu_j(A)}{A^\nu}-P^\nu_j(g)\pd{P^\mu_k(A)}{A^\nu}
\end{align*}
The Yamaguti functions $Y^\mu_{jk}$ are defined by
\begin{equation*}
6Y^\mu_{jk}(A)\=S^\mu_{jk}(A)+v^\mu_{jk}(A)+P^s_{jk}(A)
\end{equation*}
Evidently,
\begin{align*}
[S_x,S_y]
&=-x^jy^k S^\nu_{jk}(g)\pd{}{A^\nu}\\
[T_x,T_y]
&=-x^jy^k T^\nu_{jk}(g)\pd{}{A^\nu}\\
[P_x,P_y]
&=-x^jy^k P^\nu_{jk}(g)\pd{}{A^\nu}
\end{align*}
By adding the above formulae, we get
\begin{equation*}
Y(x;y)=-x^jy^k Y^\nu_{jk}(A)\pd{}{A^\nu}
\end{equation*}

\begin{lemma}
One has
\begin{subequations}
\label{stp}
\begin{align}
S^\mu_{jk}
&\=2Y^\mu_{jk}+\frac{1}{3}C^s_{jk}(S^\mu_s+2T^\mu_s)\\
T^\mu_{jk}
&\=2Y^\mu_{jk}-\frac{1}{3}C^s_{jk}(2S^\mu_s+T^\mu_s)\\
P^\mu_{jk}
&\=2Y^\mu_{jk}+\frac{1}{3}C^s_{jk}(S^\mu_s-T^\mu_s)
\end{align}
\end{subequations}
\end{lemma}

\begin{proof}
To see (\ref{uvw}a,b) use (\ref{lr-yam}a,c) . To see (\ref{uvw}c) calculate by using (\ref{lr-yam}):
\begin{align*}
[P_x,P_y]
&=[S_x+T_x,S_y+T_y]\\
&=[S_x,S_y]+[S_x,T_y]+[T_x,S_y]+[T_x,T_y]\\
&=[S_x,S_y]+2[S_x,T_y]+[T_x,T_y]\\
&=2Y(x;y)+\frac{1}{3}\left(S_{[x,y]}-T_{[x,y]}\right)
\end{align*}
and  (\ref{stp}c) easily follows.
\end{proof}

\section{Integrability conditions}

\begin{thm}
The integrability conditons of the GLE (\ref{gle_S}a--c) (\ref{gle_T}a--c) read, respectively,
\begin{subequations}
\label{gle2yam_ST}
\begin{align}
Y^s_{jk}(g)\pd{(S_gA)^\mu}{g^s}+Y^\nu_{jk}(A)\pd{(S_gA)^\mu}{A^\nu}&=Y^\mu_{jk}(S_gA)\\
Y^s_{jk}(g)\pd{(T_gA)^\mu}{g^s}+Y^\nu_{jk}(A)\pd{(T_gA)^\mu}{A^\nu}&=Y^\mu_{jk}(T_gA)
\end{align}
\end{subequations}
\end{thm}

\begin{proof}
We differentiate the GLE and use
\begin{subequations}
\label{int_ST}
\begin{align}
\frac{\p^{2}(S_gA)^\mu}{\p g^{j}\p g^{k}}&=\frac{\p^{2}(S_gA)^\mu}{\p g^{k}\p g^{j}},\quad
\frac{\p^{2}(S_sA)^\mu}{\p g^{j}\p A^{\nu}}=\frac{\p^{2}(S_gA)^\mu}{\p A^{\nu}\p g^{j}},\quad
\frac{\p^{2}(S_gA)^\mu}{\p A^{\la}\p A^{\nu}}=\frac{\p^{2}(S_gA)^\mu}{\p A^{\nu}\p A^{\la}}\\
\frac{\p^{2}(T_gA)^\mu}{\p g^{j}\p g^{k}}&=\frac{\p^{2}(T_gA)^\mu}{\p g^{k}\p g^{j}},\quad
\frac{\p^{2}(T_sA)^\mu}{\p g^{j}\p A^{\nu}}=\frac{\p^{2}(T_gA)^\mu}{\p A^{\nu}\p g^{j}},\quad
\frac{\p^{2}(T_gA)^\mu}{\p A^{\la}\p A^{\nu}}=\frac{\p^{2}(T_gA)^\mu}{\p A^{\nu}\p A^{\la}}
\end{align}
\end{subequations}
First differentiate (\ref{gle_S}a) with respect to $g^p$ and $A^\la$:
\begin{subequations}
\label{int1_S}
\begin{align}
\pd{v^s_j(g)}{g^p}\pd{(S_gA)^\mu}{g^s}
+v^s_j(g)\frac{\p^2(S_gA)^\mu}{\p g^p \p g^s}
+P^\nu_j(A)\frac{\p^2(S_gA)^\mu}{\p g^p \p A^\nu}
+\pd{T^\mu_j(S_gA)}{(S_gA)^\nu}\pd{(S_gAh)^\nu}{g^p}
&=0\\
v^s_j(g)\frac{\p^2(S_gA)^\mu}{\p A^\la \p g^s}
+\pd{P^\nu_j(A)}{A^\la}\pd{(S_gA)^\mu}{A^\nu}
+P^\nu_j(A)\frac{\p^2(S_gA)^\mu}{\p Ah^\la \p A^\nu}
+\pd{T^\mu_j(S_gA)}{(S_gA)^\nu}\pd{(S_gA)^\nu}{A^\la}
&=0
\end{align}
\end{subequations}
Now multiply (\ref{int1_S}a) by $w^p_k(g)$ and (\ref{int1_S}b) by $u^p_k(g)$ and add the resulting formulae. On the right hand side of the resulting formula use again the GLE (\ref{gle_S}a); then transpose the indexes $j$ and $k$ and subtract the result from the previous one. Then it turns out that due to (\ref{int_ST}a) all terms with the second order partial derivatives vanish and result reads
\begin{equation}
\label{gle2a_S}
v^{s}_{jk}(g)\pd{(S_gA)^{\mu}}{g^{s}}+P^{\nu}_{jk}(A)\pd{(S_gA)^{\mu}}{A^{\nu}}=T^{\mu}_{jk}(S_gA)
\end{equation}
By acting analogously with GLE (\ref{gle_S}b,c) we get
\begin{subequations}
\label{gle2b_S}
\begin{align}
u^{s}_{jk}(g)\pd{(S_gA)^{\mu}}{g^{s}}+T^{\nu}_{jk}(A)\pd{(S_gA)^{\mu}}{A^{\nu}}=P^{\mu}_{jk}(S_gA)\\
w^{s}_{jk}(g)\pd{(S_gA)^{\mu}}{g^{s}}+S^{\nu}_{jk}(A)\pd{(S_gA)^{\mu}}{A^{\nu}}=S^{\mu}_{jk}(S_gA)
\end{align}
\end{subequations}
Now add (\ref{gle2a_S}), (\ref{gle2b_S}a) and (\ref{gle2b_S}b) to obtain (\ref{gle2yam_ST}a).
 
It remains to show that  (\ref{gle2a_S}), (\ref{gle2b_S}a) and (\ref{gle2b_S}b)  are equivalent to (\ref{gle2yam_ST}a).
By using (\ref{stp}a--c) calculate
\begin{align*}
v^{s}_{jk}(g)\pd{(S_gA)^{\mu}}{g^{s}}+P^{\nu}_{jk}(A)\pd{(S_gA)^{\mu}}{A^{\nu}}&-T^{\mu}_{jk}(S_gA)
\overset{(\ref{uvw}\text{b}), (\ref{stp}\text{b,c})}{=}\\
u^{s}_{jk}(g)\pd{(S_gA)^{\mu}}{g^{s}}+T^{\nu}_{jk}(A)\pd{(S_gA)^{\mu}}{A^{\nu}}&-P^{\mu}_{jk}(S_gA)
\overset{(\ref{uvw}\text{a}), (\ref{stp}\text{b,c})}{=}\\
w^{s}_{jk}(g)\pd{(S_gA)^{\mu}}{g^{s}}+S^{\nu}_{jk}(A)\pd{(S_gA)^{\mu}}{A^{\nu}}&-S^{\mu}_{jk}(S_gA)
\overset{(\ref{uvw}\text{c}), (\ref{stp}\text{a})}{=}\\
&=2\left(Y^s_{jk}(g)\pd{(S_gA)^\mu}{g^s}+Y^\nu_{jk}(A)\pd{(S_gA)^\mu}{A^\nu}-Y^\mu_{jk}(S_gA)\right)
\end{align*}
Integrability conditions (\ref{gle2yam_ST}b) are proved analogously by using (\ref{int_ST}b).
\end{proof}

\section*{Acknowledgement}

Research was in part supported by the Estonian Science Foundation, Grant 6912.

\bigskip\noindent
Department of Mathematics\\
Tallinn University of Technology\\
Ehitajate tee 5, 19086 Tallinn, Estonia\\ 
E-mail: eugen.paal@ttu.ee

\end{document}